\documentclass[12pt]{article}
\usepackage{amsmath,amstext,amsbsy,amssymb}
\newcommand{\fr}{\frac}
\newcommand{\lb}{\label}
\newcommand{\ti}{\tilde}
\newcommand{\fb}{\begin{equation}}
\newcommand{\son}{\end{equation}}
\newcommand{\ab}{\begin{align}}
\newcommand{\as}{\end{align}}
\newcommand{\abs}{\begin{aligned}}
\newcommand{\ass}{\end{aligned}}
\newcommand{\ba}{\begin{array}}
\newcommand{\ea}{\end{array}}

\newcommand{\la}{\lambda}

\newcommand{\ep}{\epsilon}

\newcommand{\Ac}{{\cal A}}
\newcommand{\Lc}{{\cal L}}

\newcommand{\Cc}{{\cal C}}
\newcommand{\Rc}{{\cal R}}
\newcommand{\Bc}{{\cal B}}
\newcommand{\Rcc}{\mathbb{R}}
\newcommand{\Ccc}{\mathbb{C}}
\newcommand{\Ncc}{\mathbb{N}}
\newcommand{\Zcc}{\mathbb{Z}}

\newcommand{\Ic}{\mathcal{I}}

\newcommand{\Ec}{{\cal E}}

\newcommand{\dc}{{\cal D}}

\newcommand{\kd}{\delta}
\newcommand{\CD}{\Delta}
\newcommand{\epl}{\eta_+}
\newcommand{\emi}{\eta_-}
\newcommand{\ot}{\otimes}
\newcommand{\di}{\diamond}
\begin{document}
\title{}
\author{}
\date{}

\vspace{2cm}
\begin{flushright}
FGI-98-11 \\

math.QA/9809052
\end{flushright}

\vspace{1cm}
\noindent
{\Large \bf
$\mathbf{SL_q(2, \Rcc )}$ at roots of unity }

\vspace{1cm}
\noindent
{\footnotesize Hagi AHMEDOV and \"{O}mer F. DAYI}\footnote{E-mails:
 hagi@gursey.gov.tr and dayi@gursey.gov.tr.}

\noindent
{\footnotesize \it Feza G\"{u}rsey Institute,}

\noindent
{\footnotesize \it P.O.Box 6, 81220 \c{C}engelk\"{o}y--Istanbul, Turkey. }

\vspace{2cm}
\noindent
{\footnotesize {\bf Abstract.}

The quantum group $SL_q(2, \Rcc )$ at roots of unity is introduced
by means of duality pairings with the quantum algebra
$U_q(sl(2, \Rcc )).$ Its irreducible representations are constructed
through the universal $T$--matrix. An invariant integral 
on this quantum group is given. Endowed with that
some properties like unitarity and orthogonality of the irreducible
representations are discussed. }

\vspace{1.5cm}

\noindent
{\bf 1. Introduction}

\vspace{.5cm}

\noindent
One of the most interesting features of the 
quantum algebra $U_q(sl(2))$ which does not possess classical
analog is the finite dimensional cyclic irreducible
representation which appears when $q$ is a root of
unity\cite{rep}, \cite{bl}. Indeed, cyclic representations appear in
different physical applications like generalized Potts
model\cite{pot} and in classification of quantum Hall
effect wave functions\cite{dayi}.

A geometric understanding of this feature is lacking due to
the fact that structure of the related quantum group
$SL_q(2)$ at roots of unity is not well established,
although, there are  encouraging results in this
direction\cite{dun}--\cite{spe}.
When $q$ is not a root of unity $SL_q(2)$ and $U_q(sl(2))$
are duals of each other\cite{frt},\cite{vil}.
Hence, it would be reasonable
to extend this property  to obtain $SL_q(2)$
when $q$ is a root of unity. However, this is not
straightforward, because when $q$ is a root of unity
the usual dual brackets become to be ill defined. 
To cure this shortcoming
one should alter the usual number of variables 
taking part in the duality relations. Then one can  
define the quantum group $SL_q(2)$  at roots of unity.
Obviously, this fact should be reflected in $U_q(sl(2))$
at roots of unity such that the number of variables 
needed to define it should be changed consistently.

Our aim is to clarify  construction of $SL_q(2)$
at roots of unity as dual of $U_q(sl(2))$
and study them in terms
of the usual representation theory techniques. Because of the
involutions adopted, indeed we work with
$SL_q(2, \Rcc )$ and $U_q(sl(2,\Rcc )).$

In the sequel we  first discuss
$U_q(sl(2,\Rcc ))$ and $SL_q(2, \Rcc )$
for a generic $q$  
in terms of some new variables which
are suitable to define orthogonal
duality pairings. Then, we  discuss degeneracies
arising in dual brackets when we
deal with $q^p=1$ for an odd integer $p$
and present a procedure for  getting
rid of them. This yields the definition of
$SL_q(2, \Rcc )$ at roots of unity, whose subgroups are
also studied.

Once the concepts are clarified we first study irreducible
representations of 
$U_q(sl(2,\Rcc ))$ 
and then work out
the universal \mbox{$T$--matrix}. These representations as well as the
\mbox{$T$--matrix} are utilized to find out irreducible representations
of  $SL_q(2, \Rcc )$ at roots of unity.
Finally, we give the
definition of invariant integral on
$SL_q(2, \Rcc )$ and discuss the related structure
of the representations like unitarity and orthogonality.

\vspace{.5cm}
\noindent
{\bf 2. $\mathbf{U_q(sl(2, \Rcc ))}$ and $\mathbf{SL_q(2, \Rcc )}$
for a generic  $\mathbf{q}$}
\renewcommand{\theequation}{2.\arabic{equation}}
\setcounter{equation}{0}

\vspace{.5cm}

\noindent
The quantum algebra  $ U_q(sl(2, \Rcc )) $ is the $*$--Hopf algebra 
generated by $E_\pm$ and $K^{\pm 1}$ which satisfy the 
commutation relations
\begin{equation}\lb{u1}
KE_\pm K^{-1} =q^{\pm 1}E_\pm ,\  [E_+,E_-]=\fr{K^2-K^{-2}}{q-q^{-1}},
\end{equation}
the comultiplications
\begin{align}
\CD (E_\pm) & =  E_\pm \ot K + K^{-1} \ot E_\pm , &
\CD (K ) & = K \ot  K, &
\end{align}
the counits, the antipodes 
\begin{align}
\ep (K)& =1, &  \ep (E_\pm )& =0, &  \\
S(K)&=K^{-1},&  S(E_\pm ) & =-q^{\mp 1 }E_\pm  & 
\end{align}
and the involutions
\begin{equation}\lb{in1}
E_\pm^*=E_\pm,\ K^*=K.
\end{equation}

The quantum group $SL_q(2, \Rcc )$ is the $*$--Hopf algebra 
$A(SL_q(2, \Rcc))$ generated by  $x$, $y$, $u$ and  $v$ 
satisfying the commutation  relations
\begin{eqnarray}
ux=qxu, & vx=qxv, & yu=quy, \nonumber \\
yv=qvy, & uv=vu,  & yx-quv=xy-q^{-1}uv=1 , \lb{ncr}
\end{eqnarray}
the comultiplications 
\begin{eqnarray}
\CD x = x \ot x + u\ot v & , &
\CD u = x \ot u + u\ot y, \nonumber \\
\CD v = v \ot x + y\ot v     & , &       
\CD y = v \ot u + y\ot y, 
\end{eqnarray}
the counits, the antipodes
\begin{align}
\ep (x) & =1, & \ep (y)& =1,&  \ep (u)& =0,& \ep (v) & =0,& \\
S(x) & =y, &   S(y) & =x,&   S(u) & =-qu,&   S(v) & =-q^{-1}v& \lb{es}
\end{align}
and the involutions 
\begin{equation}\lb{in2}
x^*=x,\ y^*=y,\ u^*=u,\ v^*=v.
\end{equation}
The involutions adopted (\ref{in2}) and
the Hopf algebra operations (\ref{ncr})--(\ref{es}) imply   $|q|=1$.

Assume that there exists a
$*$--representation of $A(SL_q(2, \Rcc))$
such that $x$ admits the inverse $x^{-1}$ and the equality
\begin{equation}\lb{iden}
(1_A+q^{-1}uv)^{-1}=\sum_{k=0}^\infty (-1)^k (q^{-1}uv)^k
\end{equation}
holds.
$1_A$ and $1_U$ indicate the unit elements
of the related Hopf algebras.

Then, introduce the  new variables
\begin{equation}
\eta_+ =q^{-1/2} ux;\   \eta_- =q^{1/2} vx^{-1};\  \delta=x^2,
\end{equation}
dictated by the Gauss decomposition  
\begin{equation}
\lb{ga}
\left(
\ba{cc}
x & u \\
v & y
\ea
\right)=
\left(
\ba{cc}
1_A & 0 \\
q^{-1/2}\emi & 1_A
\ea
\right)
\left(
\ba{cc}
1_A & q^{1/2}\epl \\
0 & 1_A
\ea
\right) 
\left(
\ba{cc}
\delta^{1/2} & 0 \\
0 & \delta^{-1/2}
\ea
\right),
\end{equation}
satisfying the commutation relations
\begin{align}
\lb{hb}
\emi \epl & =q^2\epl \emi  , &
\eta_\pm \kd & =q^2 \kd \eta_\pm  .&
\end{align}
The involutions (\ref{in2}) yield
\begin{equation}
\eta_\pm^*=\eta_\pm ; \ \kd^* =\kd .
\end{equation}
Through the equality (\ref{iden}) we can define the
following Hopf
algebra operations on these variables:
\begin{align}
\CD \kd & =  \kd \otimes \kd + q^{-2}\kd^{-1}\eta_+^2\otimes
\eta_-^2\kd +(1_A+q^{-2})\eta_+\otimes \eta_- \kd, & \lb{cp1} \\
\CD \eta_+& =  \epl \ot 1_A +\kd \ot \epl + (1_A+q^2)\epl \ot \epl \emi
+q^{-2}\kd^{-1}\epl^2 \ot (1_A+q^2 \epl \emi )\emi, \lb{cp2} & \\
\CD \eta_-& =  \emi \ot 1_A +\kd^{-1} \ot \emi +
\sum_{k=1}^{\infty} 
(-1)^k q^{-k(k+1)} \kd^{-k-1} \epl^k \ot \emi^{k+1}, & \lb{cp3} \\
S(\kd ) & =  \kd^{-1} (1_A+q^{-2} \epl \emi )(1_A+\epl \emi ),\ \ \  
S(\eta_\pm )  =  -\kd^{\mp 1} \eta_\pm ,&    \\
\ep (\kd ) & =1, \ \ \  \ep (\eta_\pm )  =  0. &    \lb{ant}
\end{align}

When  $q$ is not a root of unity duality relations between
$U_q(sl(2, \Rcc ))$ and $A(SL_q(2, \Rcc ))$ are given by 
\begin{align}
\langle K^i \ ,\   \kd^{j} \rangle &=  q^{ij },&
i,j &\in \Zcc ,& \lb{dbna} \\
\langle E_\pm^n\ ,\  \eta_\pm^m \rangle &
= i^n q^{\pm n/2}[n]! \kd_{n,m}, & n,m & \in \Ncc ,  \lb{dbn}
\end{align}
where 
\[
[n]=\fr{q^n-q^{-n}}{q-q^{-1}}
\]
is the $q$--number.

\vspace{.5cm}
\noindent
{\bf 3. $\mathbf{U_q(sl(2, \Rcc ))}$ and $\mathbf{SL_q(2, \Rcc )}$
when $\mathbf{q^p=1}$}
\renewcommand{\theequation}{3.\arabic{equation}}
\setcounter{equation}{0}

\vspace{.5cm}
\noindent
When $q^p=1 $ (we deal with $p=$odd integer) 
for any integer $j$
we have the conditions $q^{jp}=1$ and $[jp]=0.$
So that, the dual brackets (\ref{dbna}) and  (\ref{dbn})
are degenerate. To remove the degeneracy in (\ref{dbna}) 
we put the restrictions
\begin{align}
\lb{rdb}
K^p & =1_U, & \kd^p &=1_A.&
\end{align}
By means of these conditions and the new variables
\[
\dc (m) \equiv \fr{1}{p}\sum_{l=0}^{p-1}q^{-lm}\kd^l ,
\]
instead of (\ref{dbna}) we  have 
\begin{equation}\lb{dbnb}
\langle K^n \ ,\   \dc (m) \rangle =  \kd_{n,m},\ n,m\in[0, p-1].
\end{equation}

Removing the degeneracies in (\ref{dbn}) can be achieved  in
terms of the 
following two procedures.
Take $m,n \in [0,p-1] $ in (\ref{dbn}). Let
\noindent
\begin{equation}
\lb{pgr}
\eta_\pm^p=0 ,
\end{equation}
but introduce the new variables
\begin{equation}
\lb{zpm}
z_\pm\equiv\lim_{q^p=1} \fr{\eta^p_\pm}{[p]!},
\end{equation}
without any condition on $E_\pm .$ In the second procedure
there is no condition on $\eta_\pm$ but 
on the generators
of $U_q(sl(2,\Rcc )):$
$E_\pm^p=0$ with the new variables
$
 Z_\pm \equiv\lim_{q^p=1} \fr{E^p_\pm}{[p]!} .
$
Existence of these limits $z_\pm$ and $Z_\pm$
is discussed in (\cite{bl},
\cite{cpb} and references therein).

Although, there is one more way of defining new variables
by setting both
$E_\pm^p=0$ and $\eta_\pm^p=0$ which is studied in \cite{spe},
we will show that it can be obtained as a special case 
in our approach.

We deal with the  restrictions  (\ref{pgr}) and
the new variables  (\ref{zpm}).
Now, the duality relations are
\begin{equation}
\langle E_\pm^n\ ,\  \eta_\pm^m \rangle = i^n q^{\pm n/2}[n]! \kd_{n,m},
\ \ n,m\in[0,p-1]
\end{equation}
and 
\begin{align}
\langle {\cal E}_\pm^s\ ,\  z_\pm^t \rangle
& = i^s s! \kd_{s,t}, & s,t\in \Ncc ,
\end{align}
where $\Ec_\pm \equiv (-1)^{\fr{p+1}{2}}E^p_\pm .$
Obviously, $z_\pm $  commute with the
other elements and satisfy the Hopf algebra operations
\begin{align}
S(z_\pm ) & =  - z_\pm ,& \ep (z_\pm )& =0, &
 z_\pm^* & =z_\pm , &
\end{align}
\begin{align}
\CD z_+ & = z_+ \ot 1_A +1_A \ot z_+ +
\sum_{k=1}^{p-1}\fr{q^{k^2}}{[k]![p-k]!}
\epl^{p-k} \kd^{k} \ot (-q^2 \epl\emi ; q^{2})_{(p-k)}  \epl^k, &  
\lb{cp4}\\
\CD z_- & =  z_- \ot 1_A +1_A \ot z_- +
\sum_{k=1}^{p-1}\fr{q^{-k^2}}{[k]![p-k]!}
\emi^{p-k} \kd^{-k} (-\epl\emi ; q^{-2})_k \ot \emi^k,  \lb{cp5}  &
\end{align}
where we used the notation
\[
(a;q)_k \equiv \prod_{j=1}^k (1-aq^{j-1}) .
\]
Let,
$SL_q(2,\Rcc |p)$ denotes the $*$--Hopf algebra
$A(SL_q(2,\Rcc |p))$ generated by  $\eta_\pm$
and $\delta$ through the Hopf structure given
by (\ref{hb})--(\ref{ant}). Due to the restrictions
(\ref{rdb}) and (\ref{pgr})
$SL_q(2,\Rcc |p)$ is a finite group with dimension
$p^3$.

When we deal with any $f(z_+,z_-)\equiv f(z)\in
C^\infty (\Rcc^2)$ (the space of all infinitely
differentiable functions on $\Rcc^2$)
\begin{equation}
\lb{cpf}
\CD (f(z))=f(z_0)+f_+^\prime (z_0) c_+
+f_-^\prime (z_0)c_-+f_{+-}^{\prime \prime } (z_0) c_+c_-,
\end{equation}
where $z_0\equiv (z_+\ot 1_A+1_A\ot z_+, z_-\ot 1_A+1_A\ot z_- )$ and
$c_\pm $ are given by the remaining terms of (\ref{cp4}), (\ref{cp5})
which are nilpotent $c_\pm^2=0.$
Here, $f_\pm^\prime (z_0)$
and $f_{+-}^{\prime \prime } (z_0)$
indicate derivatives of $f$ with respect to $z_\pm$ and $z_+z_-$
evaluated at $z_0.$
We also have
\begin{eqnarray}
\lb{otf}
S(f(z))=f(-z), & \ep (f(z))=f(0), & f(z)^*=\overline{f(z)},
\end{eqnarray}
where bar indicates complex conjugation.

We are ready to give the definition:

\noindent
{\it $SL_q(2,\Rcc )$
at roots of unity $(q^p=1)$ is the $*$--algebra
$A (SL_q(2, \Rcc ))= A (SL_q(2,\Rcc |p))
 \times C^\infty (\Rcc^2)$
possessing the Hopf algebra structure given by
(\ref{hb})--(\ref{ant}) and (\ref{cpf})--(\ref{otf}). }

Let  the convolution product
$\xi :A\rightarrow V$ be a homomorphic map of the Hopf algebra $A$ 
onto the linear space $V$. We set 
\begin{eqnarray}
\lb{di}
\xi\diamond g = (id \otimes \xi )\Delta (g), & g\di \xi =
(\xi\otimes id )\Delta (g), &
\xi\di \xi = (\xi\otimes \xi )\Delta . 
\end{eqnarray}
$\xi\di g$ and $g\di \xi$ belong to $A\otimes V$ and $V\otimes A,$
respectively; $\xi\di \xi$ is a homomorphic map of
$A\otimes A$ onto $V\otimes V.$

Obviously, $SL_q(2,\Rcc |p)$ is an invariant subgroup of 
$SL_q(2,\Rcc )$ at roots of unity. Moreover, in terms of the 
homomorphisim $\xi_c : A(SL_q(2,\Rcc ))\rightarrow C^\infty(\Rcc^2):$
\begin{eqnarray}
\xi_c (\eta _{\pm }) = 0, & \xi_c (\kd )=1, &  \xi_c (z_\pm )=z_\pm ,
\end{eqnarray}
one can observe that the comultiplication (\ref{cpf}) yields
\begin{equation}
\lb{xic}
\xi_c \di \xi_c (f (z) )=f(z_0).
\end{equation}
Written on the coordinates $z_\pm :$
\begin{equation}
\label{tra}
\xi_c \di \xi_c (z_\pm )= z_{\pm }\otimes 1_A+1_A\otimes z_{\pm } ,
\end{equation}
indicates that $*$--Hopf algebra $C^\infty(\Rcc^2)$  
is the translation group which is a subgroup 
of the $SL_q(2,\Rcc )$ at roots of unity.

There is another subgroup $SO(1,1 | p)$,
given in terms of the homomorphism
\begin{equation}
\xi_t (\eta _{\pm })=0,\ \ \ \xi_t (\kd )=t,
\end{equation}
where $t^p=1.$ The right sided coset 
$\Cc^{(1,1)}_q=SL_q(2,\Rcc | p) /SO(1,1| p)$
is the subspace of $A(SL_q(2,\Rcc | p))$ defined by
\begin{equation}
A(\Cc^{(1,1)}_q)= \{ g\in A(SL_q(2, \Rcc |p)): \ \ \  \xi_t \di g 
=g\otimes 1_A\}. 
\end{equation}
One can show that
\begin{equation}
\xi_t \di \eta _{+}^n\eta _{-}^m\kd ^k = 
\eta _{+}^n\eta _{-}^m\kd^k\otimes t^k. 
\end{equation}
So that,
$\eta _{+}^n\eta _{-}^m$,  $n$,
$m\in [0, p-1]$ form a basis of $A(\Cc_q^{(1,1)})$. 
Observe that
\begin{equation}
\lb{iba}
e^\pm_{nm}=\frac{\eta^{p-1-n}_+\eta^{p-1-m}_-
\pm \eta^n_+\eta^m_- }{\sqrt{q^{2n+1}+q^{-2n-1}}}, \ \ \ \ \ n,m\in [0,p-1]
\end{equation}
define  a basis which are independent
in the range
\begin{equation}
\label{dom}
n\in [0, n_0-1], \ \ \ m\in [0,2n_0 ]; \ \ \ \ \ n=n_0, \ \ \ m\in [0,n_0 ], 
\end{equation}
where $n_0=\fr{p-1}{2}.$
The number of independent elements of
$e^+_{nm}$ and $e^-_{nm}$ are
$\frac{p^2+1}{2}$ and $\frac{p^2-1}{2}$.
The quantum hyperboloid 
$H_q^{(1,1)}=SL_q(2, \Rcc)/SO(1,1|p)$ 
is defined through
the subspace of $A(SL_q(2, \Rcc ))$
\begin{equation}
\lb{qhy}
A(H_q^{(1,1)}) =A(\Cc^{(1,1)}_q) \times C^\infty (\Rcc^2) .
\end{equation}

The homomorphism
\begin{align}
\xi_l (\epl )&=\eta, & \xi_l (\emi ) & =0,& \xi_l (\kd ) & = t, &
\end{align}
defines another subgroup of $SL_q(2,\Rcc )$
denoted by $E_q(1).$ Its Hopf algebra structure
is inherited from that of   $A(SL_q(2,\Rcc )).$
The right sided coset
$ \Rcc_q  =SL_q(2,\Rcc |p )/E_q(1) $  is 
given through  the subspace
\begin{equation}
\lb{nss}
A(\Rcc_q )= \{ g\in A(SL_q(2|p)): \ \ \  \xi_l \di g =g\otimes 1_A\}. 
\end{equation}
Observe that elements of this space are polynomials in $\emi .$

We should also define:

\noindent
{\it The quantum algebra
$U_q(sl(2,\Rcc ))$ at roots of unity is generated by $E_\pm$,
$\Ec_\pm$ and $K$ with the restriction
$K^p=1_U.$ Its basis elements are 
\begin{align*}
& \Ec_+^s \Ec_-^tE^m_+E_-^nK^k & n,m,k& \in [0,p-1],& s,t &\in \Ncc .
\end{align*}
Its $*$--Hopf algebra structure is given by (\ref{u1})--(\ref{in1})
and }
\begin{align*}
\CD (\Ec_\pm ) & =\Ec_\pm \ot 1_U + 1_U\ot \Ec_\pm, &
S(\Ec_\pm ) & =-\Ec_\pm, & \ep (\Ec_\pm ) & =0, & 
\Ec_\pm^* & = \Ec_\pm .&
\end{align*}

In terms of the homomorphism $\xi_a :U_q(sl(2,\Rcc ))\rightarrow
U_q(sl(2,\Rcc |p )) $
\begin{align*}
\xi_a (E_\pm ) & = E_\pm , & \xi_a (K) &=K, & 
\xi_a (\Ec_\pm ) &=0, &
\end{align*}
we can define $U_q(sl(2,\Rcc |p))$
the sub--Hopf algebra of $U_q(sl(2,\Rcc ))$
generated by
\begin{align*}
E_\pm^p & =0, & K^p=1_U.
\end{align*}
Obviously, the discrete quantum algebra
$U_q(sl(2,\Rcc |p))$ is in non--degenerate duality with
$SL_q(2,\Rcc |p)$.
This is the case  studied in \cite{spe}.

\vspace{.5cm}

\noindent
{\bf 4. Irreducible $\mathbf{*}$-representations of
$\mathbf{U_q(sl(2,\Rcc ))}$
when $\mathbf{q^p=1}$}
\renewcommand{\theequation}{4.\arabic{equation}}
\setcounter{equation}{0}

\vspace{.5cm}

\noindent
The homomorphism $\Lc^\la:\ U_q(sl(2))\rightarrow {\rm Lin}\ A(SO(1,1|p))$ 
given by 
\fb
\begin{aligned}
\Lc^\la (K)  & t^i  = q^{-i}t^i, && i  \in [0,p-1],   \\
\Lc^\la (E_-) & t^i = t^{i+1}, && i  =0,1,\cdots ,p-2,  \\
\Lc^\la (E_-) & t^{p-1} =\la_+t^0 ,&&   \\
\Lc^\la (E_+) & t^i  =  M_it^{i-1} && i  =1,\cdots ,p-1,  \\
\Lc^\la (E_+) & t^0  =  a t^{p-1} && , 
\end{aligned}
\end{equation}
where the constants are
$$
\la_-=a\prod_{i=1}^{p-1} M_i,\   M_i=a\la_+-[i-1][i],
$$
defines the cyclic irreducible representation of
$U_q(sl(2))$ ($\Bc $ type
representation)\cite{rep},\cite{bl}.

We would like to find out when  
$\Lc^\la$ defines a $*$--representation.
To this aim introduce the Hermitian form
\begin{equation}
\label{hers}
(a,b)_t=\Ic_t (a^* b ), 
\end{equation}
for $a,b \in A(SO(1,1| p))$
and the linear functional on it
\begin{equation}
\Ic_t (t^m )=\delta_{m, 0( {\rm mod}\ p)}.
\end{equation}
Moreover, we see that
\begin{align*}
e_m^{\pm } & =\frac {1}{\sqrt{2}} (t^m \pm t^{p-m} ), & m 
& \in[0,\frac{p-1}{2}] , 
\end{align*}
are orthogonal with respect to the Hermitian form (\ref{hers}): 
\begin{align*}
(e_m^\pm ,e_k^\pm )_t & = \pm \delta_{mk}, &
(e_m^\mp ,e_k^\pm )_t& =0. 
\end{align*}
Thus, with  the Hermitian form (\ref{hers})
$*$--Hopf algebra $A(SO(1,1| p))$ 
is the pseudo--Euclidean space
possessing $\frac{p+1}{2}$ positive
and $\frac{p-1}{2}$
negative signatures.

Adjoint of a linear operator is defined through
$$
(\Lc^\la (\phi )a,b)_t= (a, (\Lc^\la(\phi ) )^* b)_t,
$$
where $\phi \in U_q(sl(2,\Rcc )).$ Hence, we conclude that
if $\la_\pm$ are real $\Lc^\la$ defines a 
\mbox{$*$--representation:}
$$
(\Lc^\la (\phi ))^*= \Lc^\la(\phi^* ).
$$

The linear map  $T^{(l)}$: 
$A(\Rcc_q) \rightarrow  A(SL_q(2,\Rcc )) \times
A(\Rcc_q) $ given by
\[
T^{(l)} g(\emi )=\left( id \ot \kd^{-l}\right) \CD (\kd^lg(\emi )),
\]
for $l\in[0,\fr{p-1}{2}]$ defines irreducible representations of 
$SL_q(2, \Rcc )$.
Infinitesimal form of this global representation is 
\[
\Rc^{(l)} (\phi )g(\emi )=(\phi \ot id ) T^{(l)}g(\emi ),
\]
where $\phi \in U_q(sl(2,\Rcc )).$
We see that
\begin{align*}
\Rc^{(l)} (E_+) \emi^{l-m}& =iq^{l+1/2} [l+m]  \emi^{l-m+1}, & \\
\Rc^{(l)} (E_-) \emi^{l-m}&  =iq^{-l-1/2} [l-m]  \emi^{l-m-1}, & \\
\Rc^{(l)} (K) \emi^{l-m}&= q^m\emi^{l-m},&
\end{align*}
where $m\in [-l,l].$ These are  non--cyclic representations of
$U_q(sl(2,\Rcc ))$ ($\Ac$ type representations).

\vspace{.5cm}

\noindent
{\bf 5.  The universal T--matrix and irreducible representations of
$\mathbf{SL_q(2,\Rcc )}$ at roots of unity}
\renewcommand{\theequation}{5.\arabic{equation}}
\setcounter{equation}{0}

\vspace{.5cm}
\noindent
Let  the basis elements of 
the Hopf algebras $U(g)$ and $A(G),$
respectively, $V_a$ and $v^a$
lead to the dual brackets 
$\langle V_a,v^b \rangle=\kd_a^b,$ 
which are non--degenerate. Then the
universal $T$--matrix   $T\in U(g)\ot A(G)$
can be constructed as\cite{fg},\cite{bon}
\[
T=\sum_a V_a \ot v^a .
\]
As far as the universal
$T$--matrix is known, one can construct 
corepresentations of $A(G)$ utilizing
representations of $U(g).$

A straightforward calculation leads to the duality brackets
\begin{align}
\langle
\Ec_+^t\Ec_-^s
E_+^n  E_-^mK^k\ ,\ z_+^{t'}z_-^{s'}\eta_+^{n^\prime} \eta_-^{m'}
\dc (k')\rangle = &
i^{s+t+n+m} q^{\fr{(n-m)}{2}-nm}s!t![m]![n]!  \nonumber \\
  & \kd_{n,n'}\kd_{m,m'}\kd_{s,s'}\kd_{t,t'}\kd_{k+n+m,k'} , \lb{dbr}
\end{align}
where $n,m \in [0 , p-1].$
Therefore, the universal T--matrix can be written as
\begin{equation}
\label{utm}
T=e^{-i \Ec_+\ot z_+ -i\Ec_-\ot z_- }
\sum_{n,m,k=0}^{p-1}
\fr{i^{-n-m} q^{\fr{m-n}{2}+nm}}{[n]![m]!}
E_+^n E_-^m K^k \ot \epl^n \eta_-^m  \dc (k+n+m)  .
\end{equation}
Arranging the elements and using the cut off q--exponentials
\[
e_\pm^x=\sum_{r=1}^{p-1} \fr{q^{\pm r(r-1)/2}}{[r]!} x^r ,
\]
the universal $T$--matrix can also be written as
\begin{equation}
\label{utmm}
T=e^{-i \Ec_+\ot z_+ -i\Ec_-\ot z_- }
e_+^{i\ep_+ \ot \epl } e_-^{i\ep_-\ot \emi} D(K, \kd ),
\end{equation}
where we introduced
\begin{align*}
\ep_\pm &=-q^{\pm 1/2} E_\pm K^{-1}, & \\
D(K,\kd ) & =\fr{1}{p} \sum_{k,l =0}^{p-1} q^{-ml}K^k\ot \kd^l . &
\end{align*}
Using the explicit form (\ref{utmm}) one can show that
\begin{equation}
\label{utp}
[(*\otimes *)T]\cdot T=1_A \otimes 1_U,\ \ \ T\cdot (*\otimes
*)T=1_A \otimes 1_U .
\end{equation}
In general, $T$--matrix also satisfies 
\begin{equation}
\label{uto}
(id\otimes \Delta )T=(T\otimes 1_A)(id\otimes \sigma )(T\otimes 1_A),
\end{equation}
where $\sigma (F\otimes G)=G\otimes F$, $F$, $G\in A(SL_q(2, \Rcc )),$
is the permutation operator.

Let us illustrate how one obtains 
irreducible representations of $SL_q(2,\Rcc )$
by making use of  
the universal $T$--matrix (\ref{utm}). 
Let  $T^{(\lambda )}: A(SO(1,1| p)) \rightarrow A(SO(1,1|
p))\otimes A(SL_q(2,\Rcc ))$, 
be 
\begin{equation}
\label{repg}
\begin{aligned}
T^{(\la )}a = & e^{-i \Lc^\la (\Ec_+)\ot z_+ -i\Lc^\la (\Ec_-)\ot z_- }
e_+^{i\Lc^\la (\ep_+) \ot \epl }
e_-^{i\Lc^\la (\ep_-)\ot \emi} D(\Lc^\la (K), \kd ).
\end{aligned}
\end{equation}
Because of  (\ref{uto})  and the irreducibility of 
the representation $\Lc^\la$ we
conclude that $T^{(\la )}a$ gives a $p$--dimensional 
irreducible representation of the quantum group 
$SL_q(2,\Rcc )$
in the linear space $A(SO(1,1| p))$. 
Let us extend the Hermitian form 
(\ref{hers}) to 
\begin{equation}
\label{herss}
\{ a\ot F, b\ot G \}_t =(a,b)_t F^*G  , 
\end{equation}
where $F$, $G\in A(SL_q(2,\Rcc ))$ and $a$, $b\in A(SO(1,1| p))$. When
$\lambda_\pm$ are real numbers  the condition
(\ref{utp})  yields
\begin{equation}
\label{pseudo}
\{ T^{(\la )} a ,T^{(\la )} b \}_t = (a,b)_t 1_A. 
\end{equation}
Thus the irreducible representation 
$T^{(\la )}$ is pseudo--unitary when 
$\lambda_\pm$ are real.

We can obtain  matrix elements of the
irreducible
pseudo--unitary representations as  
\begin{equation}
D_{mn}^\la = \{  t^{p-m} \ot 1_A, T^{(\la )} t^n \}_t. 
\end{equation}
For some specific values of $n,\ m$ we performed the
explicit calculations:
\begin{equation}
\lb{00}
D_{00}^\la= e^{-i\la_+z_+-i\la_-z_-}\{ 1+\sum_{m=1}^{p-1}
\fr{(-1)^m}{ ([m]!)^2 }
\left( \prod_{j=1}^{m}M_j\right)  \rho^m \} ,
\end{equation}
where $\rho =q\epl \emi .$ For $i \neq 0,$ we obtain
\begin{eqnarray}
\lb{oi}
D_{i0}^\la & = & e^{-i\la_+z_+-i\la_-z_-} \{    
\sum_{m=0}^{p-i-1}\fr{(-1)^m i^{-i}q^{i(m-1/2)} }{[m]![m+i]!}
\left( \prod_{j=1}^{m+i}M_j  \right) \rho^m \emi^i \nonumber \\
  & &  +\sum_{m=0}^{i-1}\fr{(-1)^m i^{i-p}q^{i(p-1)/2-im} }{ [m]![p+m-i]! }
\left( \prod_{j=0}^{m}M_j  \right) \epl^{p-i}\rho^m    \}, \lb{0i}
\end{eqnarray}
where the definition $M_0\equiv \la_+$ is adopted.

The pseudo--unitarity condition  (\ref{pseudo}) implies 
\begin{equation}
(D^\la_{0 m})^*D^\la_{0 n} + \sum_{k=1}^{p-1}(D^\la_{k
m})^*D^\la_{p-k n}=(t^m, t^n)_t 1_A. 
\end{equation}
Special cases are 
\begin{eqnarray}
(D^\lambda_{00})^*D^\lambda_{00} +
\sum_{k=1}^{p-1}(D^\lambda_{0k})^*D^\lambda_{0p-k }=1_A, \nonumber \\
(D^\lambda_{0i})^*D^\lambda_{0p-i} +
\sum_{k=1}^{p-1}(D^\lambda_{k i})^*D^\lambda_{p-k p-i}=1_A.
\end{eqnarray}
Moreover, we have the addition theorem 
$$
\Delta (D^\la_{n m})=\sum_{k=0}^{p-1}D^\la_{nk} \otimes D^\la_{k m}. 
$$

\vspace{.5cm}
\noindent
{\bf 6. Regular  representation of $\mathbf{SL_q(2, \Rcc )}$ }
\renewcommand{\theequation}{6.\arabic{equation}}
\setcounter{equation}{0}

\vspace{.5cm}

\noindent
The comultiplication 
\begin{equation}
\lb{rr}
\Delta : A( SL_q(2, \Rcc ))\rightarrow A(SL_q(2, \Rcc )) \ot
A(SL_q(2, \Rcc )) 
\end{equation}
defines the  regular representation of $SL_q(2, \Rcc )$ 
in the linear space   $A\left( SL_q(2, \Rcc )\right).$
The  right and left representations of  $U_q(sl(2, \Rcc ))$ 
corresponding to the regular representation (\ref{rr})
are given, respectively, by
$$
\Rc (\phi )F\equiv \hat{\phi} F=  F\di \phi   
$$
and 
$$
\Lc (\phi )F\equiv\ti{\phi} F=  \phi \di F, 
$$
where $F\in A\left( SL_q(2, \Rcc )\right).$
Straightforward calculations yield 
the right representations
\begin{align*}
\hat{E}_+\epl^n & = 
iq^{1/2}[n]\epl^{n-1} + iq^{1/2-n}[2n]\emi \epl^{n},&
\hat{E}_+\emi^n  & =   -iq^{1/2}[n]\emi^{n+1},  \\
\hat{E}_-\emi^n & =  iq^{-1/2} [n]\emi^{n-1},  &  
\hat{K}\eta_\pm^n & =   q^{\pm n} \eta_\pm^n,  \\
\hat{E}_-\epl^n & =    0,  & \hat{E}_-\kd^n & =  0, \\
\hat{E}_+\kd^n & = i(q^{-3/2-n}+q^{-3n-7/2})[n+1]
\emi\kd^n (1-\kd_{n,0}), &
\hat{K} \kd^n &  =  q^{n} \kd^n,  \\  
\hat{E}_\pm f(z_+,z_-) & =  \fr{iq^{\pm 1/2}}{[p-1]!}\eta_\pm^{p-1}
\fr{df(z_+,z_-)}{dz_\pm }, & 
\hat{K}z_\pm & =  z_\pm ,
\end{align*}
and the left representations 
\begin{align*}
{\ti{E}}_+\epl^n &  =   iq^{n-3/2} [n]\kd \epl^{n-1}, &
{\ti{E}}_-\epl^n  & =     iq^{-n-1/2}\kd^{-1}\epl^{n+1},  \\
{\ti{E}}_-\emi^n &  =   iq^{3/2-n} [n]\kd^{-1}\emi^{n-1}, &
{\ti{K}}\eta_\pm^n  & =   \eta_\pm^n,  \\
{\ti{E}}_+\emi^n  & =   0, &
{\ti{E}}_+\kd^n  & =   0, \\
{\ti{E}}_-\kd^n  & =   iq^{3/2-n}[2n]\epl \kd^{n-1} (\kd_{n,0}-1),  &
{\ti{K}} \kd^n  & =   q^{n} \kd^n,  \\
\ti{E}_\pm f(z_+,z_-) & =  \fr{iq^{\mp 1}}{[p-1]!}\eta_\pm^{p-1}\kd^\pm
\fr{df(z_+,z_-)}{dz_\pm }, & 
\ti{K}z_\pm & =  z_\pm .
\end{align*}

Right representation of any element $\phi\in U_q(sl(2,\Rcc ))$
can be found through the above relations and making use of
the properties
\begin{align*}
\Rc (\phi \phi^\prime ) & = \Rc (\phi^\prime )\Rc ( \phi  ), \\
\hat{E}_\pm (XY) &= \hat{E}_\pm X \hat{K} Y +
\hat{K}^{-1} X\hat{E}_\pm Y,  \\
\hat{K} XY & =\hat{K}X \hat{K}Y.
\end{align*}
For left representations similar properties hold.

Although, the quantum algebra  $U_q(sl(2, \Rcc ))$ at roots of
unity possesses three Casimir elements ${\cal E}_\pm$ and 
$$
C= E_-E_+ +\fr{(qK-q^{-1}K^{-1})^2}{(q^2-q^{-2})^2},
$$
only two of them are independent. 
Thus, irreducible representations of   $U_q(sl(2, \Rcc ))$
at roots of unity are labeled by two indices.
A method of 
constructing the irreducible representations
of $U_q(sl(2, \Rcc ))$ at roots of unity
is to  diagonalize the complete set of commuting
operators   $\hat{\Ec }_\pm$, $\hat{C}$ and $\hat{K}$ on the quantum 
hyperboloid. 
Indeed, the matrices  (\ref{00}) and (\ref{0i})
can be shown to satisfy
\begin{align*}
\hat{C}D_{i0} & =a \lambda_+ D_{i0},& i & \in [0,p-1], \\
\hat{\Ec}_\pm D_{i0} &
=(-1)^{\fr{p+1}{2}}\la_\pm D_{i0}, & i & \in [0,p-1],  \\
\hat{E}_{+}D_{i0}^\la & = D_{(i+1)0}^\la ,&   i&\in [0,p-2], \\
\hat{E}_{+}D_{(p-1)0}^\la &= \la_+D_{0,0}^\la , &  &\\
\hat{E}_-D_{i0}^\la & = M_i D_{(i-1)0}^\la ,&  i &\in [1,p-1],  \\
\hat{E}_-D_{00}^\la & = aD_{0(p-1)}^\la . &  &
\end{align*}
Similar constructions can also be done in terms of
the left representations.

\vspace{.5cm}

\noindent
{\bf 7. Invariant integral on $\mathbf{ SL_q(2,\Rcc )}$
at roots of unity}
\renewcommand{\theequation}{7.\arabic{equation}}
\setcounter{equation}{0}

\vspace{.5cm}

\noindent
Recall that the invariant integral $\Ic$ on the quantum group $G_q$ 
is a linear functional  on the Hopf algebra $A(G_q)$
which for any element $a\in A(G_q)$
satisfies the left 
\begin{equation}
\lb{li} 
\Ic \di a  = 1_A \Ic (a)
\end{equation}  
and the right
\begin{equation}\lb{li1}   
a \di \Ic  = 1_A\Ic (a) 
\end{equation}
invariance conditions.

The linear functional $\Ic_p$ on the Hopf algebra 
$A(SL_q(2, \Rcc|p))$ given by 
\begin{equation}
\lb{ind}
\Ic_p  (\epl^n \emi^m \kd^k)
=q^{-1}\kd_{n,p-1}\kd_{m,p-1}\kd_{k,0 ({\rm mod} \ p)} 
\end{equation}
defines the invariant integral on the quantum group $SL_q(2, \Rcc|p)$.    
To prove that in fact the conditions (\ref{li}) and  (\ref{li1}) are
satisfied, we proceed as follows. Since $A(SL_q(2, \Rcc |p))$ is
a finite
Hopf algebra it is sufficient 
to show that
(\ref{li}) and  (\ref{li1})
are satisfied after taking their dual pairings:
\begin{eqnarray}
\Ic_p  \left(\Rc (\phi )P\right) & = & \Ic_p  (P) \ep (\phi ), \lb{rp1}\\
\Ic_p  \left(\Lc (\phi )P\right) & = & \Ic_p  (P) \ep (\phi )  ,\lb{rp2}
\end{eqnarray}
for all elements $\phi\in U_q(sl(2, \Rcc ))$
and $P\in A(SL_q(2,\Rcc |p)).$
One can show that
\begin{align}
\Ic_p  (\hat{E}_\pm \epl^n \emi^m \kd^k ) & =  0, &  
\Ic_p  (\hat{K} \epl^n \emi^m \kd^k ) & = 
\Ic_p  ( \epl^n \emi^m \kd^k ), &\lb{rpm}  \\
\Ic_p  ({\ti{E}_\pm} \epl^n \emi^m \kd^k ) & =  0, & 
\Ic_p  ({\ti{K}} \epl^n \emi^m \kd^k ) & =
\Ic_p  ( \epl^n \emi^m \kd^k ). &\lb{lk} 
\end{align}
Moreover, for any two elements  $\phi_1,\ \phi_2$ right and left
representation satisfy the relations
\begin{eqnarray*}
\Ic_p  ({\cal R} (\phi_1 \phi_2) P) & = & \ep (\phi_1 \phi_2) \Ic_p  (P), \\
\Ic_p  (\Lc (\phi_1 \phi_2) P) & =  & \ep (\phi_1 \phi_2) \Ic_p  (P).
\end{eqnarray*}
Therefore, (\ref{rp1}) and (\ref{rp2}) are satisfied. This 
leads to the conclusion that
(\ref{ind}) is the invariant integral on $SL_q(2, \Rcc |p).$

Observe that 
\begin{equation}
\Ic_p ( P^* )= \overline{\Ic_p (P)}
\end{equation}
and define the Hermitian form $(\cdot ,\cdot )_p$
on the quantum group $SL_q(2,\Rcc | p)$ as
\begin{equation}
\label{phm}
(P,Q)_p= \Ic _p (PQ^* ). 
\end{equation}
The basis elements $e^\pm_{nm}$ (\ref{iba}) of
$A(\Cc_q^{(1,1)})$ are orthonormal in terms the above form:  
$$
(e^\pm_{nm},e^\pm_{n^\prime m^\prime})_p=\pm
\delta_{nn^\prime}\delta_{mm^\prime}, \ \ \ (e^\pm_{nm},e^\mp_{n^\prime
m^\prime})_p=0. 
$$
Any element $\pi \in A (\Cc^{(1,1)}_q)$ can be represented as
\begin{equation}
\pi =\sum_{nm} \pi^+_{nm}e^+_{nm} + \sum_{nm} \pi^-_{nm}e^-_{nm}, 
\end{equation}
where $\pi^\pm_{nm}\in \Ccc$ and $n$, $m$ take values in 
the domain (\ref{dom}). Then, the pseudo--norm of $\pi$
\begin{equation}
(\pi ,\pi )_p=\sum_{nm} \pi^+_{nm}\overline {\pi^+_{nm}} - \sum_{nm}
\pi^-_{nm} \overline{\pi^-_{nm}},
\end{equation}
shows that the metric of the space $A(\Cc^{(1,1)}_q)$ possesses
$\frac{p^2+1}{2}$ positive and $\frac{p^2-1}{2}$ negative signatures.

We should also define invariant integral on the translation subgroup
for being able to obtain it on $SL_q(2,\Rcc ).$

Let $C_0^\infty (\Rcc^2)$ be the space of all
infinitely differentiable functions with
finite support in $\Rcc^2$.
The linear functional on  $C_0^\infty(\Rcc^2):$
\begin{equation}
\lb{incr}
\Ic_c (f)= \iint_{-\infty}^\infty dz_+ dz_- f (z_+,z_- ).
\end{equation}
where $f\in C_0^\infty (\Rcc^2),$
is clearly the invariant integral on the translation group
satisfying 
\begin{align}
\label{inv}
(\Ic_c\otimes id)(\xi_c\di \xi_c)(f)& = \Ic_c  (f), &
(id \otimes \Ic_c)(\xi_c\di \xi_c)(f) &= \Ic_c  (f). &
\end{align}

Let $A_0(SL_q(2,\Rcc ))$  be the subspace of $A (SL_q(2,\Rcc ))$
defined as 
\begin{equation}
A_0(SL_q(2,\Rcc )) =  C_0^\infty(\Rcc^2)\times A(SL_q(2,\Rcc | p)) 
\end{equation}
and  $\Ic_w$ be the linear functional acting on it as
\begin{equation}
\label{haar}
\Ic_w(F)=\sum_n\Ic_p(P_n)\Ic_c(f_n),
\end{equation}
where $F=\sum_nP_nf_n$ and $f_n\in C_0^\infty (\Rcc^2),\
P_n\in A (SL_q(2,\Rcc | p))$. Let us prove that $\Ic_w$ is
the invariant integral on  $A_0(SL_q(2,\Rcc )).$
On an element  $G=Pf$ we have 
\begin{equation}
\lb{ns1}
\Ic_w\di G  =  (id \otimes \Ic_w)\CD (P)\CD (f).
\end{equation}
One can observe from (\ref{xic}) that any function $f(z)$
evaluated at $z=z_0$ can be written as
\[
f(z)|_{z_0}=\xi_c\di \xi_c (f(z)) .
\]
Hence, (\ref{ns1}) yields
\[
\Ic_w\di G =(id\otimes \Ic_w)\left\{ \CD (P) \left[
\xi_c\di \xi_c (f)+
 c_+ \xi_c\di \xi_c (f_+^\prime )
+c_- \xi_c\di \xi_c (f_-^\prime )
+ c_+c_- \xi_c\di \xi_c (f_{+-}^{\prime \prime } ) \right] \right\}.
\]
by making use of (\ref{cpf}). Moreover, the properties of the
invariant integrals (\ref{ind}), (\ref{incr}) and
(\ref{haar}) permits us to write
\begin{equation}
\Ic_w\di G  =  \Ic_p(P)\Ic_c (f)+
(id\otimes \Ic_p)\left\{ \Delta (P) \left[
 c_+ \Ic_c (f_+^\prime ) 
 +c_- \Ic_c (f_-^\prime )
+ c_+c_-, \Ic_c (f_{+-}^{\prime \prime } )\right] \right\}.
\end{equation}
Because $f\in C_0^\infty (\Rcc^2),$ we have
\begin{equation}
\Ic_c(\frac {df}{dz_{\pm }})=
\Ic_c(\frac{d^2f}{dz_{+}dz_{-}})=0.
\end{equation}
Hence,
\begin{equation}
\Ic_w\di G=\Ic_p(P)\Ic_c(f)=\Ic_w(G),
\end{equation}
which together with the linearity of the functional
$\Ic_w$ implies
\begin{equation}
\Ic_w\di F=\Ic_w(F)\ \ 
\text{for any}\ \ F\in A_0(SL_q(2,\Rcc )). 
\end{equation}         
Right invariance condition can  be proved similarly.
Therefore,
$\Ic_w$ is the invariant integral on the quantum group
$SL_q(2,\Rcc )$ at roots of unity.

Let us introduce the bilinear form
\begin{equation}
\label{her}
(F,G)_w=\Ic_w(FG^*),
\end{equation}
where $F,G \in A_0(SL_q(2,\Rcc ))$, which is Hermitian because
\begin{equation}
\Ic_w(F^{*})=\overline{\Ic_w(F)}.
\end{equation}
Consider the subspace of $A_0(SL_q(2,\Rcc ))$
\begin{equation}
A_0(H_q^{(1,1)}) = C_0^\infty (\Rcc^2)\times A(\Cc_q^{(1,1)}),
\end{equation}
whose arbitrary element $X$ can be written as
\begin{equation}
X=\sum_{nm}f_{nm}^{+}e_{nm}^{+}+
\sum_{nm}f_{nm}^{-}e_{nm}^{-},
\end{equation}
where $e_{nm}^{\pm}$ are given by
(\ref{iba}) in the domain (\ref{dom}).
We then have
\begin{equation}
(X,X)_w=\sum_{nm}
\Ic_c(f_{nm}^{+}\overline{f_{nm}^{+}})-
\sum_{nm}\Ic_c(f_{nm}^{-}\overline{f_{nm}^{-}}).
\end{equation}
Thus, $A_0(H_q^{(1,1)})$ endowed with the Hermitian form (\ref{her})
is a pseudo--Euclidean space.

The  comultiplication
\[
\CD :A_0 (H_q^{(1,1)} )\rightarrow A_0 (SL_q(2,\Rcc )) 
\ot A_0(H_q^{(1,1)} ),
\]
defines the left quasi--regular representation of $SL_q(2, \Rcc )$ 
in  $A_0(H_q^{(1,1)} )$.
Let us extend the Hermitian form $(\ , \ )_w$ to 
$\{\ , \ \}_w$ by setting
\[
\{ F\ot X, G \ot Y \}_w \equiv F G^*(X,Y)_w,
\]
where $F,\ G \in A_0(SL_q(2,\Rcc ))$ and 
$X,\ Y \in A_0(H_q^{(1,1)}).$   
We have
\begin{equation}
\lb{hog}
\{\CD(X), \CD(Y)\}_w =1_A (X,Y)_w,
\end{equation}
which implies that the  left quasi--regular representation is 
pseudo--unitary.

For any  $\phi\in U_q(sl(2, \Rcc ))$ and $F\in A_0(SL_q(2,\Rcc ))$ 
the duality brackets satisfy the property
\[
\overline{\langle \phi^*,F\rangle } = \langle \phi, (S(F))^*\rangle,
\]
which together with  the pseudo--unitarity condition (\ref{hog}) implies
\[
(\Rc (\phi )X,Y)_w= (X,\Rc (\phi^* )Y)_w.
\]
Thus, the antihomomorphism $\Rc: U_q(sl(2, \Rcc ))\rightarrow {\rm Lin}
A_0(H^{(1,1)}_q)$ given in Section 6 defines the $*$--representation of 
the quantum algebra in the pseudo--Euclidean space 
$A_0(H^{(1,1)}_q)$.

Note that the matrix elements of the pseudo--unitary irreducible
representations (\ref{00}), (\ref{oi})
satisfy the orthogonality condition 
\[
(D_{n0}^\la, D_{m0}^{\la^\prime})_w=\kd (\la_+-\la_+^\prime )
\kd (\la_--\la_-^\prime ) N_i\kd_{n+m, 0({\rm mod} \ p)},
\]
where  $N_n$ are some  normalization constants.


\newpage

\begin{thebibliography}{99}
\bibitem{rep}G. Lusztig, {\it Adv. Math.} {\bf 70 } (1988) 237;\\
M. Rosso, {\it Commun. Math. Phys.} {\bf 117 } (1988) 581;\\
P. Roche and D. Arnaudon, {\it Lett. Math. Phys.} {\bf 17 } (1989) 295;\\
C. de Concini and  V. G. Kac, {\it Representations of Quantum 
Groups at Root of Unity,} Progress in Mathematics Vol. 92 
(Birkh\"{a}use, Boston, 1990).

\bibitem{bl}L. C. Biedenharn and M. A. Lohe, {\it Quantum
Group Symmetry and q--Tensor Algebras,}
(World Scientific, Singapore, 1995).

\bibitem{pot}E. Date , M. Jimbo, K. Miki and T. Miwa,
{\it Commun. Math. Phys. }{\bf 137} (1991) 133; \\
V.V. Bazhanov, R. M. Kashaev, V. V. Mangazeev and Yu. G. Stroganov,
{\it Commun. Math. Phys. }{\bf 138} (1991) 393.

\bibitem{dayi}\"{O}. F. Dayi,
{\it J.  Phys. A: Math. Gen.} {\bf 31}  (1998) 3523.

\bibitem{dun}R. S. Dunne, {\it Intrinsic anyonic spin through
deformed geometry}, DAMTP/97-19, hep-th/9703137.

\bibitem{spe}D. V.  Glushenkov and A. V. Lyakhovskaya,
{\it Regular
representation of the quantum Heisenberg double  $\{ U_q(sl(2)),\  
Fun(SL(2))\}$ (q is a root of unity ),}  UUITP-27/1993,
hep-th/9311075,\\
R. Coquereaux, A. O. Garcia and R.Trinchero,
{\it Differential Calculus and Connections on a Quantum Plane at a Cubic 
Root of Unity,} CPT-98/P.3632, math-ph/9807012. 

\bibitem{frt}L. D. Faddeev, N. Yu. Reshetikhin and L. A. Takhtajan,
{\it Leningrad Math. Journal } {\bf 1} (1989) 178, \\
T. Masuda, K. Mimachi, Y. Nagakami, M, Noumi and K. Ueno,
{\it Lett. Math. Phys. }, {\bf 19} (1990) 187.

\bibitem{vil}N. Ja. Vilenkin and A. O. Klimyk, 
{\it Representation of Lie Groups and Special Functions, vol 3},
(Dordrecht: Kluwer Akad. Publ. 1992).

\bibitem{cpb}V.  Chari and A. Pressley, {\it  Quantum Groups}
(Cambridge: Camb. Univ. Press, 1994).

\bibitem{fg}C. Fronsdal, A. Galindo, 
{\it Lett. Math. Phys. } {\bf 27} (1993) 59.

\bibitem{bon}F. Bonechi, E. Celeghini, R. Giachetti,
C.M. Pere\~{n}a, E. Sorace and M. Tarlini, 
{\it J. Phys. A: Math. Gen.} {\bf 27} (1994) 1307; \\
F. Bonechi, N. Ciccoli, R. Giachetti,
E. Sorace and M. Tarlini, 
{\it Commun. Math. Phys.} {\bf 175} (1996) 161.


\end{thebibliography}
\end{document}